\documentclass[preprint,12pt]{elsarticle}
\usepackage{graphics}
\usepackage{url}
\usepackage{multirow}
\graphicspath{{./}{figure/}}

\usepackage{amssymb}
\usepackage{amsthm}
\usepackage{amsfonts}
\usepackage{amsmath}
\usepackage{amssymb}

\newcommand{\RR}{\field{R}}
\newcommand{\NN}{\field{N}}

\newtheorem{theorem}{\bf Theorem}[section]

\newtheorem{proposition}[theorem]{\bf Proposition}

\newcommand{\Proofend}{\hfill$\diamondsuit$}

\def\NN{{\mathbb N}}
\def\RR{{\mathbb R}}

\def\I{{\mathcal{I}}}

\newcommand{\PSNR}{\rm{PSNR}}
\newcommand{\MSE}{\rm{MSE}}
\newcommand{\SSIM}{\rm{SSIM}}
\newcommand{\cov}{\rm{cov}}

\journal{ArXiv}

\begin{document}
\begin{frontmatter}
\title{{L}agrange--{C}hebyshev Interpolation for image resizing}
\tnotetext[mytitlenote]{The research has been accomplished within Research ITalian network on Approximation (RITA) and Approximation Theory research group  of Unione Matematica Italiana (TA--UMI). It has been partially supported by
GNCS--INdAM and University of Basilicata (local funds).
}
\author[address1, address]{Donatella Occorsio}
\ead{donatella.occorsio@unibas.it}
\author[address]{Giuliana Ramella}
\ead{giuliana.ramella@cnr.it}
\author[address]{Woula Themistoclakis\corref{cor1}}
\ead{woula.themistoclakis@cnr.it}
\cortext[cor1]{Corresponding author}
\address[address1]{Department of Mathematics and Computer Science, University of Basilicata, Viale dell'Ateneo Lucano 10, 85100 Potenza, Italy}
\address[address]{C.N.R. National
Research Council of Italy, Institute for Applied Computing ``Mauro Picone'', Via P. Castellino, 111, 80131 Naples, Italy}
\begin{abstract}
 Image resizing is a basic tool in image processing and in literature we have many methods, based on different approaches, which are often specialized in only upscaling or downscaling. 
In this paper, independently of the (reduced or enhanced) size we aim to get,
we approach the problem at a continuous scale where the underlying continuous image is globally approximated  by the tensor product Lagrange polynomial interpolating at a suitable grid of first kind Chebyshev zeros. This is a well--known approximation tool that is widely used in many applicative fields, due to the optimal behavior of the related Lebesgue constants. Here we show how Lagrange--Chebyshev interpolation can be fruitfully applied also for resizing an arbitrary digital image in both downscaling and upscaling.
The performance of the proposed 
 method has been tested in terms of the standard SSIM and PSNR metrics. The results indicate that, in upscaling, it is almost comparable with the classical Bicubic resizing method with slightly better metrics, but in downscaling a much higher performance has been observed in comparison with Bicubic and other recent methods too. Moreover, in downscaling cases with an odd scale factor, we give an estimate of the mean squared error produced by our method and prove it is theoretically null (hence PSNR equals to infinite and SSIM equals to one) in absence of noise or initial artifacts on the input image.
\end{abstract}
\begin{keyword}
Image resizing\sep Image downscaling\sep Image upscaling \sep Lagrange interpolation \sep Chebyshev nodes
\MSC  94A08\sep 68U10\sep 41A05 \sep 62H35
\end{keyword}
\end{frontmatter}
\section{Introduction}\label{sec:intro}
In this paper, we deal with the problem of image resizing. This has been widely investigated over the past decades and is still an active research area
characterized by many applications in different domains, including image transmission, satellite image
analysis, gaming, remote sensing, etc. (see e.g. \cite{Atkinson,quattro,cinque,Meijering,venticinque,sei})
In literature, downscaling and upscaling are often considered  separate problems (see e.g. \cite{diciannove,ventiquattro,venticinque}), and most of the existing methods are specialized in only one direction,  sometimes for a limited range of scaling factors (see e.g. \cite{scaling1,ventitre}).  The method we are going to introduce works in both down and up scaling directions and for "large" scaling factors as well. It falls into the class of the interpolation methods and it is based on a non-standard modeling of the image resizing problem.

In order to introduce the adopted model we premise that, from the mathematical viewpoint, a continuous image $\I$ is a function $f(x,y)$ of the spatial coordinates which, without losing the generality, we can assume belonging to the open square $A=]-1,1[^2$. Hence, its digital version $I$ of $n\times m$ pixels is supposed to be made of the values that $f$ takes on a discrete grid of  nodes $X_{n\times m}\subset A$.
 Regarding such grid, it is generally supposed that $X_{n\times m}=X_n\times X_m$ where we set
 \begin{equation}\label{X}
X_\mu:=\{x_k^\mu \ : \ k=1,\ldots,\mu\}\subset ]-1,1[,\qquad \forall \mu\in\NN.
\end{equation}
A typical and natural choice of the (univariate) system of nodes $X_\mu$ is to divide $[-1,1]$ in $(\mu+1)$ equal parts and to take the $\mu$ internal equidistant nodes, i.e.
\begin{equation}\label{equi-nodes}
X_\mu^{equ}=\left\{-1+\frac{2k}{\mu+1} \ :\ k=1,\ldots, \mu\right\}.
\end{equation}
However, it is well known that equally spaced  nodes are not the best choice for Lagrange interpolation since they lead to exponentially growing Lebesgue constants, whereas optimal Lebesgue constants (growing at the minimal projection rate) are provided by the Chebyshev nodes of the first kind   
\begin{equation}\label{x-Che}
x_k^\mu= \cos\left[\frac{(2k-1)\pi}{2\mu}\right],\qquad k=1,\ldots,\mu,\qquad \forall\mu\in\NN,
\end{equation}
(see e.g. \cite{mastromilo, Tref} for a short excursion on the topic).

Discrete polynomial approximation based on Chebyshev zeros  is a pillar in approximation theory and practices. It has been widely studied and applied in several fields and also recently, Chebyshev--like grids such as Xu points and Padua points have been introduced to get optimal interpolation processes on the square (see e.g. \cite{Xu, Xu1, Padua, Padua1, occothem_AMC}).  Nevertheless, to our knowledge, its usage in image processing has been mainly limited to particular cases, such as Magnetic Particle Imaging that is strictly related to Lissajous curves generating the Padua points ( see e.g. \cite{MPI, MPI1, MPI2}).

\begin{figure}[h]
\label{fig-nodes}
\begin{centering}
\includegraphics[scale=0.15]{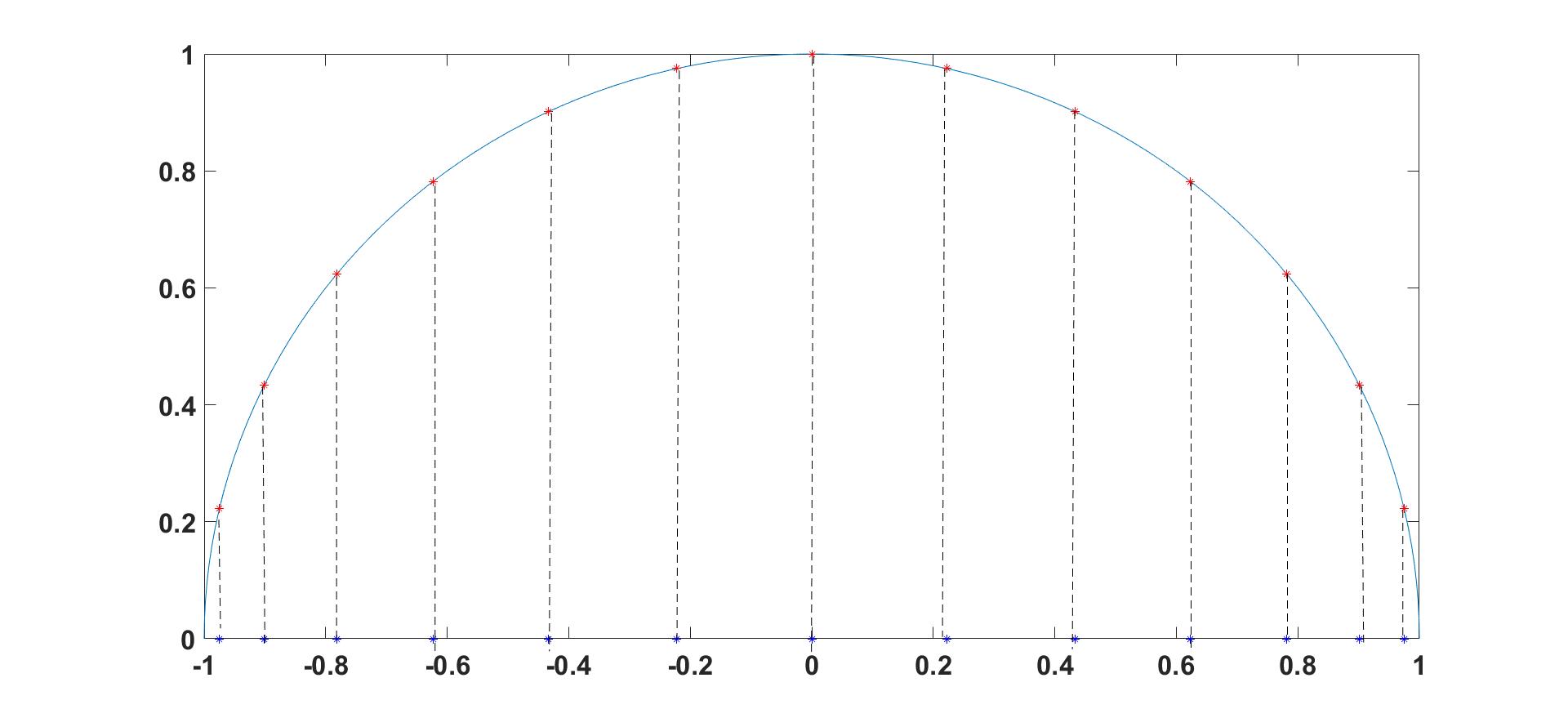}
\caption{The first kind Chebyshev zeros of order $\mu=15$ (black stars) and their arccosine (red stars)}
\end{centering}
\end{figure}
Indeed, 
at a first look, Chebyshev zeros may seem unsuitable for sampling an arbitrary image since they are not equidistant on the segment $]-1,1[$ of each spatial coordinate, but they are {\sl arc sine} distributed, becoming denser at the endpoints $\pm 1$ (see Figure \ref{fig-nodes}). Nevertheless,  in our model, we transfer the sampling question from the segment $]-1,1[$ to the semicircle of radius 1 centered at the axes origin. Thus, using the standard setting $t=\arccos x$, instead of the usual equidistant nodes on the segment, our method takes the following equidistant nodes on the semicircle
\begin{equation}\label{t-equi}
t_k^\mu:= \frac{(2k-1)\pi}{2\mu},\qquad k=1,\ldots,\mu,\qquad \mu\in\NN,
\end{equation}
These nodes define the Chebyshev zeros $x_k^\mu=\cos t_k^\mu$ in (\ref{x-Che}) and the grid
\begin{equation}\label{XxX}
X_{n\times m}=\left\{(x_i^n,x_j^m):\quad i=1:n, \ j=1:m\right\}, \qquad n,m\in\NN
\end{equation}
which we actually choose to sample arbitrary (continuous) images $\I$ and to obtain digital images of $ n \times m $ pixels, for arbitrary $n,m\in\NN$.

 According to this model,  let $I=f(X_{n\times m})$ be the digital image of $n\times m$ pixels, i.e.
 \begin{equation}\label{I}
 I_{i,j}:=f(x_i^n,x_j^m),\qquad i=1:n,\ j=1:m,
 \end{equation}
and let $I^{in}$ be its more or less corrupted version really available as input data, say
\begin{equation}\label{I-in}
I^{in}_{i,j}:=\tilde f(x_i^n,x_j^m),\qquad i=1:n,\ j=1:m,
 \end{equation}
where $\tilde f$ denotes a 
corrupted version of the continuous image. Starting from the data $I^{in}$, the  resizing problem aims to find  a good approximation of the resampled image $I^{res}$ consisting of the samples of $f$ at the different (more or less dense) grid $X_{N\times M}$, i.e.
\begin{equation}\label{I-res}
I^{res}_{k,h}:=f(x_k^N,x_h^M),\qquad k=1:N,\ h=1:M.
 \end{equation}
Hence, under our modeling, the resizing problem raises the following approximation question: starting from the approximate values of $f$ at the Chebyshev grid $X_{n\times m}$, how to approximate $f$ at another (coarser in downscaling or finer in upscaling) Chebyshev grid $X_{N\times M}$? The solution we propose   relies on the global approximation of a function by its bivariate Lagrange  polynomial interpolating at the Chebyshev grid $X_{n\times m}$.

This well--known interpolation polynomial can be easily deduced via tensor product from the univariate case (see for instance \cite{occorusso2011,occothemLNCS,apnum, mata}) and fast algorithms can be implemented for its computation ( see e.g. \cite{occothem_AMC}). Moreover, we recall that wavelets techniques can be applied to such kind of interpolation \cite{Kilgore, Prestin, ThFischer}.

Hence, in order to get the resampled image, instead of $f$, we propose to sample a suitable Lagrange interpolation polynomial at the new grid. 
More precisely, we use the pixels of the input image $I^{in}$ to build the Lagrange-Chebyshev polynomial $L_{n,m}\tilde f$ interpolating $\tilde f$ at $X_{n\times m}$ and we get the output image of the desired $N\times M$ size (here denoted by $I^{out}$) by sampling $L_{n,m}\tilde f$ at the  mesh $X_{N\times M}$.

To test the performance of such Lagrange--Chebyshev Interpolation (LCI) method we have considered the usual SSIM (Structured Similarity Index Measurement)
and PSNR (Peak Signal to Noise Ratio) metrics as computed by MatLab, and  we have carried out the experimentation   on  several kinds of images collected in 5  datasets, having different  characteristics on the  size and the quality of the images, the variety of subjects, etc.

 For an initial comparison, we focused on interpolation methods working as well in both down and up scaling directions.  In general, these methods find a resized image by local or global interpolation of the continuous image $f$, using in some way the initial $n\times m$ pixels, namely the more or less corrupted values of $f$ at the grid of nodes $X_{n\times m}$.
Typically, for such grid it is assumed the equidistant nodes model (\ref{equi-nodes}) and most of the interpolation methods act locally on each pixel of the input image. This is the case of the classical bicubic interpolation method BIC \cite{keys} implemented by the Matlab built--in  function  \texttt{imresize} which we have used
 either in upscaling and downscaling to get a first comparison.

 By the extensive experimentation we carried out, we assess that the performance  of the LCI method is superior in downscaling (d-LCI method), while in  upscaling (u-LCI method) there are only slight differences and the two methods seem almost comparable.

 Thus it  seemed relevant to us to further investigate only the downscaling case by testing d-LCI method in comparison with two other recent downscaling methods, here denoted by DPID \cite{Weber} and L$_0$  \cite{Liu}, both of them not belonging to the family of interpolation methods for downscaling. The effective better performance in downscaling has been confirmed also in these cases, with an increasing gap as the downscaling factor increases. Moreover,  by running the public available Matlab code, we have observed that  DPID and L$_0$ methods work only for integer scale factors $n/N$ and $m/M$ and they are not always implementable for images  with large sizes due to too long running times or too much memory. On the contrary, LCI method keeps the runtimes contained and offers high flexibility since it works in both downscaling and upscaling  with non integer scale factors too, being possible to specify the desired final size or the scale factor as input parameters. We remark that similar nice features are shared by the Matlab function \texttt{imresize}, but with lower performances  in downscaling. Moreover, in downscaling, for all odd scale factors, we state a theoretical estimate of the Mean Squared Error (MSE) obtained by our d-LCI method arriving to prove that it is null if $I^{in}=I$. Such a theoretical result has been confirmed by our experimentation that also shows it does not hold for \texttt{imresize} applied to the same input images.

The paper is organized as follows. For a simpler exposition of the idea behind the LCI method, we  first introduce the method for a monochrome image in Section 2 and then, in Section 3, we expose it for RGB color images.  Both  monochrome and color cases are considered in Section 4 where some advantages of our approach are discussed. Finally, Section 5 consists of all the numerical experiments.
\section{Re-sampling  gray-level images}\label{section1}
In the case of a monochrome image,  $\mathcal{I}$ is represented at a continuous scale by a scalar function $f(x,y)$ whose values are the grey levels at the spatial coordinates $(x,y)\in A$. Moreover, at a discrete scale, digital images of $n\times m$ pixels are matrices $I\in \RR^{n\times m}$ whose elements are the samples of $f$ at the discrete nodes set $X_{n\times m}$ defined in our model by (\ref{XxX}) and (\ref{x-Che}). Then, according to the notation introduced in the previous section, the resizing problem of $\mathcal{I}$ (in upscaling or in downscaling, respectively) consists in finding a good approximation of the matrix $I^{res}\in \RR^{N\times M}$ given by (\ref{I-res}) starting from the matrix (having a reduced or  enlarged size, respectively) $I^{in}\in \RR^{n\times m}$ defined as in (\ref{I-in}) being $\tilde f$ a more or less corrupted version of $f$.

We approach such an approximation problem  by using the following bivariate (tensor product) Lagrange--Chebyshev polynomial based on the available data $I^{in}$
 \begin{equation}\label{lag}
L_{n,m}\tilde f(x,y):=\sum_{i=1}^n\sum_{j=1}^m \tilde f(x_i^n, x_j^m)\ell_i^n(x) \ell_j^m(y), \qquad (x,y)\in A,
\end{equation}
where, for all $\mu\in\NN$, $\ell_k^\mu$ denotes  the $k-$th fundamental Lagrange polynomial related to the nodes system $X_\mu$, namely
 \begin{equation}\label{fund-lag}
\ell_k^\mu(\xi):=\prod_{\scriptsize\begin{array}{c}
s=1\\ [-.1cm] s\ne k\end{array}}^\mu\frac{\xi- x_s^\mu}{x_k^\mu- x_s^\mu}, \qquad \xi\in [-1,1],\qquad,\qquad k=1:\mu.
\end{equation}
It is well--known that the Lagrange--Chebyshev polynomial in (\ref{lag}) interpolates $\tilde f$ at the grid $X_{n\times m}$, i.e.
 \begin{equation}\label{inter}
L_{n,m}\tilde f(x_i^n, x_j^m)= \tilde f(x_i^n, x_j^m)=I^{in}_{i,j}, \qquad i=1:n,\ j=1:m.
\end{equation}
The samples of such polynomial at the re-scaled grid $X_{N\times M}$ constitute the approximate resized image provided by the LCI method.
Hence, the LCI method  computes  the matrix $I^{out}\in \RR^{N\times M}$ whose entries are the following
\begin{equation}\label{Iout}
I^{out}_{k,h}:=L_{n,m}\tilde f(x_k^N,x_h^M)
, \qquad k=1:N,\quad h=1:M.
\end{equation}
Introducing the Vandermonde-like matrices
\begin{equation}\label{V}
V_1:=\left[\ell_i^n(x_k^N)\right]_{i,k}\in \RR^{n\times N}, \qquad V_2:=\left[\ell_j^m(x_h^M)\right]_{j,h}\in \RR^{m\times M},
\end{equation}
by (\ref{lag}) and (\ref{Iout}), the output matrix $I^{out}$ can be computed from the input matrix $I^{in}$ according to the following matrices identity
\begin{equation}\label{L-prod}
I^{out}=V_1^T I^{in} V_2,
\end{equation}
Note that in case we have to resize a lot of images for the same fixed sizes, formula (\ref{L-prod}) also allows working in parallel with pre-computed matrices $V_i$ .

Moreover, it is well known that the fast computation of the matrices $V_i$ can be achieved by using, instead of  (\ref{fund-lag}), the following more convenient form of the fundamental Lagrange polynomial
\begin{equation}\label{fund-trig}
\ell_k^\mu (\cos t)=\frac 2\mu \sum_{r=0}^{\mu-1} { }^\prime\cos\left[\frac{(2k-1)r\pi}{2\mu}\right]
\cos\left[r t\right],\qquad k=1:\mu,\quad t\in[0,\pi],
\end{equation}
where, as usual, the prime on the summation symbol means that the first addendum is halved.
Hence, by (\ref{fund-trig}) we get that the matrices $V_1$ and $V_2$ can be computed by using fast cosine transform algorithms (see e.g. \cite{tasche}) being
\begin{eqnarray}\label{V1}
(V_1)_{i,k}=\ell_i^n(x_k^N)\hspace{-.2cm}&=&\hspace{-.2cm}\frac 2n \sum_{r=0}^{n-1} { }^\prime\cos\left[\frac{(2i-1)r\pi}{2n}\right]
\cos\left[\frac{(2k-1)r\pi}{2N}\right],\\
\nonumber
&&\qquad i=1: n,\ k=1:N,\\
\label{V2}
(V_2)_{j,h}=\ell_j^m(x_h^M)\hspace{-.2cm}&=&\hspace{-.2cm}\frac 2m \sum_{s=0}^{m-1} { }^\prime\cos\left[\frac{(2j-1)s\pi}{2m}\right]
\cos\left[\frac{(2h-1)s\pi}{2M}\right],\\
\nonumber&&\qquad j=1: m,\  h=1:M.
\end{eqnarray}
\section{Re-sampling RGB color images}\label{section2}
Now we are going to introduce the method for a general color image $\I$.
In this case, behind $\I$ there is a vector function $\mathbf{f}:A\rightarrow \RR^3$ whose components $\mathbf{f}=(f_R,f_G,f_B)$ represent $\mathcal{I}$  in the RGB color space.  According to our model, the input image $\mathbf{I}^{in}$ and the target resized image $\mathbf{I}^{res}$ are represented in the RGB space as the follows
\[
\begin{array}{ll}
\mathbf{I}^{in}\in \RR^{n\times m\times 3}, &\mathbf{I}^{in}\equiv (I^{in}_R, I^{in}_G, I^{in}_B)\\
\mathbf{I}^{res}\in \RR^{N\times M\times 3}, & \mathbf{I}^{res}\equiv (I^{res}_R, I^{res}_G, I^{res}_B),
\end{array}
\]
where we assume that
\begin{equation}\label{I_in_vector}
\begin{array}{lll}
\left(I^{in}_R\right)_{i,j}={\tilde f_R}(x_i^n,x_j^m), &\ & \left(I^{res}_R\right)_{k,h}={f_R}(x_k^N,x_h^M),\\
\left(I^{in}_G\right)_{i,j}={\tilde f_G}(x_i^n,x_j^m), &\ & \left(I^{res}_G\right)_{k,h}={f_G}(x_k^N,x_h^M), \\
\left(I^{in}_B\right)_{i,j}={\tilde f_B}(x_i^n,x_j^m), &\ & \left(I^{in}_B\right)_{k,h}={f_B}(x_k^N,x_h^M).
\end{array}
\end{equation}
holds for all $ [i,j]=[1:n,1:m]$ and $[k, h]=[1:N,1:M]$.

Thus, similarly to the monochrome case and with obvious meaning of the notation, we start from
 \begin{equation}\label{interp3}
\mathbf{I}^{in}_{i,j}=\tilde{\mathbf{f}}(x_i^n,x_j^m)
\qquad i=1:n,\quad j=1:m,
\end{equation}
and approximate $\mathbf{I}^{res}$ by
\begin{equation}\label{interp3out}
\mathbf{I}^{out}_{k,h}=L_{n,m}\tilde{\mathbf{f}}(x_k^N,x_h^M)
\qquad k=1:N,\quad \ h=1:M
\end{equation}
i.e., for $ k=1:N$ and $h=1:M $, we define
\begin{equation}\label{approx3}
\mathbf{I}^{out}_{k,h}=\left(L_{n,m}{\tilde f_R}(x_k^N,x_h^M),\ L_{n,m}{\tilde f_G}(x_k^N,x_h^M),\ L_{n,m}{\tilde f_B}(x_k^N,x_h^M)\right).
\end{equation}
Hence, by applying the same argument of the monochrome case, we get that the RGB components of the output image $\mathbf{I}^{out}\equiv (I^{out}_R, I^{out}_G, I^{out}_B)$ are given by
$$\mathbf{I}^{out}=(V_1^T I^{in}_R V_2,\ V_1^T I^{in}_G V_2,\ V_1^T I^{in}_B V_2),$$
with $V_1,V_2$ defined in (\ref{V1})--(\ref{V2}).
\section{Model analysis}
For a general analysis of the error we get in approximating the target resized image
$I^{res}=[ {\mathbf{f}} (x_i^n, x_j^m)]_{i,j}$ with the output image $I^{out}=[ \L_{n,m}\tilde{\mathbf{f}} (x_h^N, x_k^M)]_{h,k}$ produced by LCI method, we refer the reader to the wide existing literature on the Lagrange interpolation error estimates (see e.g. \cite{mastromilo} and the references therein). In this paper, we are interested to measure the performance of LCI method through some standard metrics usually used in Image Processing (see e.g. \cite{Ra-metrics} and the references therein).
 In particular, we focus on the Mean Squared Error (MSE) defined as follows
\begin{equation*}
\MSE(I^{res},I^{out})= \begin{cases} \displaystyle\frac{1}{MN}\|I^{res}-I^{out}\|_F^2, \hspace{3.5cm}\textrm{ gray-level images},\\ \\
\frac 1 {3NM} \left\{ \|I_R^{res}-I_R^{out}\|_F^2+\|I_G^{res}-I_G^{out}\|_F^2+\|I^{res}_B-I_B^{out}\|_F^2\right\}, \\ \hspace{7cm}\textrm{ RGB color images}, \end{cases}
\end{equation*}
being  $\|\cdot \|_F$  the Frobenius norm
\[
\|A\|_F:=\left(\sum_{k=1}^N\sum_{h=1}^M A_{k,h}^2\right)^\frac 12, \qquad A=(A_{k,h})\in\RR^{N\times M}.
\]
Moreover, we consider the Peak Signal to Noise Ratio (PSNR) defined by the previous MSE as follows
\begin{equation}\label{psnr}
\PSNR(I^{res},I^{out})=20 \displaystyle \log_{10}\left( \frac{\max_f}{\sqrt{\MSE(I^{res},I^{out})}}\right),
\end{equation}
with $\max_f=255$, since we use the   representation by 8 bits per sample.

Finally, in our experiments we evaluate the Structured Similarity Index Measurement (SSIM). For grey-level images it is defined   as
\begin{equation}\label{def_ssim}
\SSIM(I^{res},I^{out})=  \displaystyle
\frac{\left[2\mu(I^{res}) \mu(I^{out})+c_1\right]\left[2\cov(I^{res},\ I^{out})+c_2\right]}
{\left[\mu^2(I^{res})+\mu^2(I^{out})+c_1\right]\left[\sigma^2(I^{res})+\sigma^2(I^{out})+c_2\right]},
\end{equation}

where  $\mu(A), \sigma(A)$ and $\cov (A, B)$  indicate the average, variance and covariance, respectively, of the matrices $A,B$, and $c_1, c_2$ are constants usually fixed as $c_1=(0.01\times L), c_2=(0.03\times L)$ with the dynamic range of the pixel values $L=255$ for 8-bit images.

In the case of RGB color images, making the conversion  to the  color space  YCbCr, the $\SSIM$ is computed by the above definition applied to the intensity Y (luma) channel.


Based on the previous metrics, the results of  wide experimentation will be reported in the next section. In the sequel, we are going to underly two particular features of LCI method.

A first aspect is related to the model choice which, for odd  downscaling factors $s:=n/N=m/M>1$, leads to state the following
\begin{proposition}\label{prop}
Under the previously introduced notation, if there exists $\ell\in\NN$ s.t. $n=(2\ell-1)N$ and $m=(2\ell-1)M$ hold,  then we have
 \begin{equation}\label{dim_mse}
 \MSE(I^{res},I^{out})\le s^2 MSE(I, I^{in}),\qquad s:=(2\ell -1).
 \end{equation}
Moreover, if  $I=I^{in}$ then we get the best results
  \begin{equation}\label{dim_ssim}\MSE(I^{res},I^{out})=0,\qquad\mbox{and}\qquad \SSIM (I^{res},I^{out})=1.\end{equation}
\end{proposition}
{\it Proof. } Let's give the proof in the case of a gray-level image since for RGB color images the statement will easily follow  by considering the single RGB components.

As a keynote of the proof, we observe that for any $\ell\in\NN$, we have
\[
n=(2\ell-1)N \qquad\Longrightarrow\qquad X_N\subset X_n.
\]
More precisely, setting $s:=(2\ell-1)$, it is easy to check that
\begin{equation}\label{nodi_3}
n= s N\qquad \Longrightarrow \qquad  x_k^N=x_i^n  \quad\mbox{with}\quad i=\frac{s(2k-1)+1}2,
\quad k=1:N,
\end{equation}
where we remark that for all $k=1:N$, the numerator $s(2k-1)+1$ is certainly even since $s$ is odd.

Consequently, we deduce (\ref{dim_mse}) from  (\ref{nodi_3}) as follows
\begin{eqnarray*}
&&\MSE(I^{res},I^{out}) =\frac 1{NM}\sum_{k=1}^N\sum_{h=1}^M\left[f\left(x_k^N,x_h^M\right)-L_{n,m}\tilde f\left(x_k^N,x_h^M\right)\right]^2\\
&=&\frac 1{NM}\sum_{k=1}^N\sum_{h=1}^M\left[f\!\left(x_{\frac{s(2k-1)+1}2}^n,x_{\frac{s(2h-1)+1}2}^m\right)\!- \!
\tilde f\!\left(x_{\frac{s(2k-1)+1}2}^n,x_{\frac{s(2h-1)+1}2}^m\right)\right]^2\\
&\le&
\frac 1{NM}\sum_{i=1}^n\sum_{j=1}^m\left[f\left(x_i^n, x_j^m\right)-\tilde f\left(x_i^n, x_j^m\right)\right]^2=s^2 MSE(I,I^{in}).
\end{eqnarray*}
Moreover, in the case that $I=I^{in}$, we have that $L_{n,m}f=L_{n,m}\tilde f$ and, due to the nesting property $X_{N\times M}\subset X_{n\times m}$, by (\ref{inter}), for any $h=1:N$ and $k=1: M$, we get
\begin{equation}\label{identity}I^{out}_{h,k}=
L_{n,m}\tilde f(x_h^N,x_k^M)=L_{n,m}f(x_h^N,x_k^M)
=f(x_h^N,x_k^M)=I^{res}_{h,k},
\end{equation} that implies the first equation in (\ref{dim_ssim}). Finally, (\ref{identity}) also yields
$$2\cov(I^{res},\ I^{out})=\sigma(I^{res})^2+\sigma(I^{out})^2,$$
$$2\mu(I^{res})\mu(I^{out})=\mu(I^{res})^2 +\mu(I^{out})^2$$
which conclude  the proof of (\ref{dim_ssim}).
\Proofend

We remark that the previous proposition is a direct consequence of  our choice of the sampling model based on Chebyshev zeros
 (\ref{x-Che}), instead of the usual equally spaced nodes (\ref{equi-nodes}).

\begin{figure}[h]
\begin{center}
\includegraphics[scale=0.25]{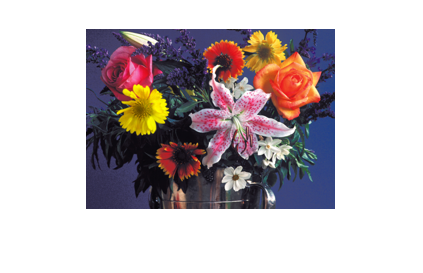}
\includegraphics[scale=0.25]{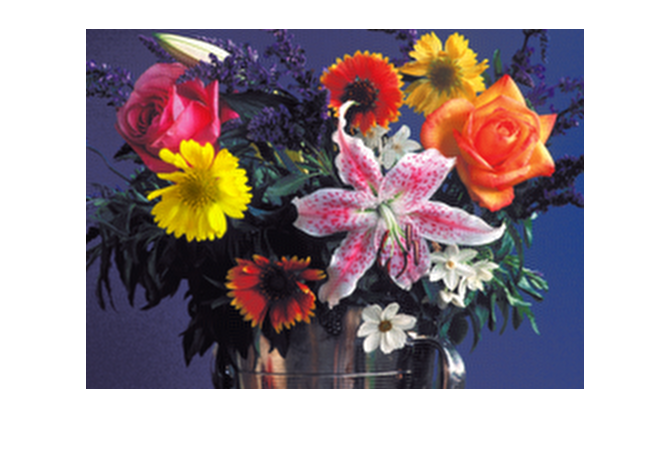}
\includegraphics[scale=0.25]{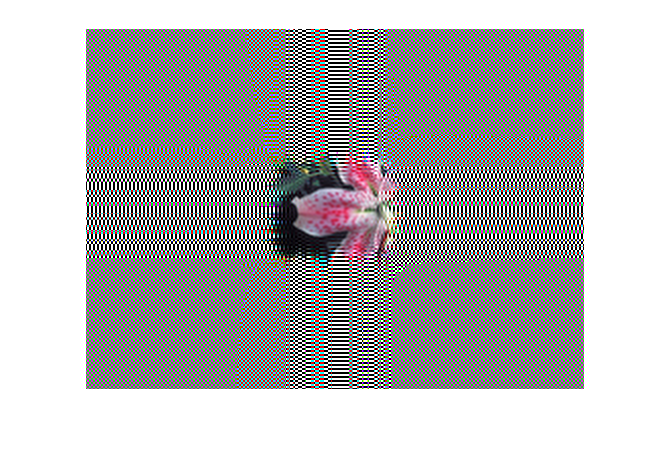}
\caption{The well--known image {\it Flowers} taken from USC--SIPI \cite{SIPI} (left) upsampled at the double scale by  Lagrange interpolation at  Chebyshev nodes (middle) and at equally spaced nodes (right)}
\label{fig:lenax2}
\vspace{-.5cm}
\end{center}
\end{figure}

 Such choice (and more generally the choice of "good" interpolation knots) allows to use global interpolation processes that are not possible in the case of equally spaced nodes.
In fact, as second aspect, we aim to highlight the relevance of  the choice of  good interpolation knots by showing  
the effects we have in the  performance of global Lagrange interpolation methods associated with the univariate nodes set $X_{\mu}^{equ}$ in (\ref{equi-nodes}).
The disastrous effects of the exponential growth of Lebesgue constants are visible in the next experiment concerning an image zooming $\times 2$ from the size $n\times n$ with $n=256$ to the size $N\times N$ with $N=512$.
Figure \ref{fig:lenax2} displays  the images obtained starting from the same input image (Fig.\ref{fig:lenax2}, left) and using the Lagrange interpolation polynomial based on  equally spaced nodes $X_n^{equ}\times X_n^{equ}$ (Fig.\ref{fig:lenax2}, right) and on the Chebyshev grid $X_{n\times n}$ (Fig.\ref{fig:lenax2}, middle). The images we see are the samples of the previous Lagrange polynomials at the respective $(N\times N)$--grid and we can check how the equally--spaced sampling model (\ref{equi-nodes})  yields wrong results outside of a small region of the output image.
\section{Experimental results}\label{section3}
In this section, we propose a selection of the  extensive experimentation  we carried on  to test  the benefits of our procedure in downscaling (d-LCI method) and in upscaling (u-LCI method), comparing our results with other resizing techniques.
\subsection{\it Comparison methods}
As already mentioned, we compare LCI with BIC  provided by MatLab \texttt{imresize} with \texttt{'bicubic'} option where we recall that a new  pixel is determined  from a
weighted average of the 16 closest pixels using an interpolating  cubic convolution function satisfying prescribed assumptions of smoothness on $f$  to gain at least a cubic order of convergence \cite{keys}. The comparison BIC--LCI regards either down and up resizing cases,  giving as input either the scale factor or the final size of the desired image. Note that we have also tested the other different options of \texttt{imresize} (namely \texttt{'linear'} and \texttt{'nearest neighborhood'}) but for brevity, we do not report the results since they give no new insight.

Moreover, to further investigate the performance in downscaling case, we  compare d-LCI  with other two recent downscaling methods not belonging to the interpolation method class. These methods are briefly indicated as  DPID method (described in  \cite{Weber} with the code available at \cite{codeDPID} ) and  L$_0$ method 
(described in  \cite{Liu}, with the code available at \cite{codeL0}).
To be fair,  we remark that we forced DPID and L$_0$ also in upscaling direction for several scaling factors (the choice of the desired size is not allowed by DPID and L$_0$) but the results were very poor w.r.t. BIC and LCI  and they have been not reported here.

Finally, we point out that all  the previous methods  have run on the same PC with  Intel Core i7 3770K CPU @350GHz configuration.
\subsection{\it  Datasets}
Since the results of any method 
depend on the type of input image, we have conducted the experimentation   on  different kinds of 8--bits images collected in five publicly available  datasets, whose characteristics   (acronym, number of images, range of their sizes) are  synthesized in  Table \ref{tab:1}.
\begin{table}[h]
\caption{Datasets list}
\label{tab:1}
\centering{\begin{tabular}{|l|c|c|}
\hline
\textbf {Dataset} & n. of Images & Sizes \\
\hline
BSDS500& 500 & $481\times 321$ or $321\times 481$ \\
NASA & 17 & from 500$\times$334 to 6394$\times$3456  \\
YAHOO & 96 & from 500$\times$334 to 6394$\times$3456 \\
13US & 13 & from 241$\times$400 to 400$\times$310   \\
URBAN100 & 100 & from 1024$\times$564 to 1024$\times$1024  \\
\hline
\end{tabular}
}
\end{table}

We point out that we have chosen the datasets in Table \ref{tab:1} since they are the same considered in  \cite{Weber,Liu} for DPID and L$_0$ methods.  In particular:
\begin{itemize}
\item
BSDS500 dataset \cite{Martin}, available at \cite{Berkeley}, has been used for testing L$_0$ in \cite{Liu}. It is  sufficiently general and provides a large variety of images often employed for testing many other methods with different image analysis tasks such  as image segmentation \cite{segmentation1,segmentation2,segmentation3,segmentation4}),   color quantization \cite{40,41,42,43}, etc.
\item
NASA Image Gallery \cite{NASA} and  YFCC100M (Yahoo Flickr Creative Commons 100 Million)  \cite{Thomee} datasets (here briefly denoted by NASA and YAHOO, respectively) are  used by DPID  in \cite{Weber}. 
\item
13US \cite{Oztireli} contains natural images, available at \cite{13US} and originally taken from the MSRA Salient Object Database \cite {Liu2011}. It is another dataset used in \cite{Weber}. 
\item
 URBAN100 dataset \cite{Urban100},
 concerning urban scenes with images having one dimension equal to 1024, is commonly used  to evaluate the performance of super-resolution models. It is used as test images for L$_0$ in \cite{Liu}.
\end{itemize}
All the previous datasets have been considered for either up and down scaling.

We remark that in all the performed experiments the input images are not available from the datasets. Hence, following a typical approach in image quantitative evaluation, we fix all the images from the datasets as target images (i.e. $I^{res}$ in our notation) and we apply  BIC to generate the input image (namely $I^{in}$) common to all methods. More precisely, to run  [$\times s$] upscaling methods ([$:s$] downscaling methods, resp.) we generate the input image $I^{in}$  by applying BIC  to $I^{res}$ in the opposite [$:s$] ([$\times s$], resp.) scaling direction.
\newline
Note that, according to such procedure, in testing [$\times s$] upscaling methods  we may find that $I^{res}$ does not have both the dimensions $N$ and $M$ that are divisible for $s$. In this case, in order to generate $I^{in}$ we use \texttt{imresize} with the size option, requiring $n=\lfloor \frac N s\rfloor$ and  $m=\lfloor \frac M s\rfloor$. Moreover, once obtained $I^{in}$, both BIC and u-LCI  run by specifying again the size $N\times M$ instead of the scaling factor.
\subsection{\it Visual comparison}
\begin{figure}[h]
\begin{center}

\includegraphics[width=10 cm]{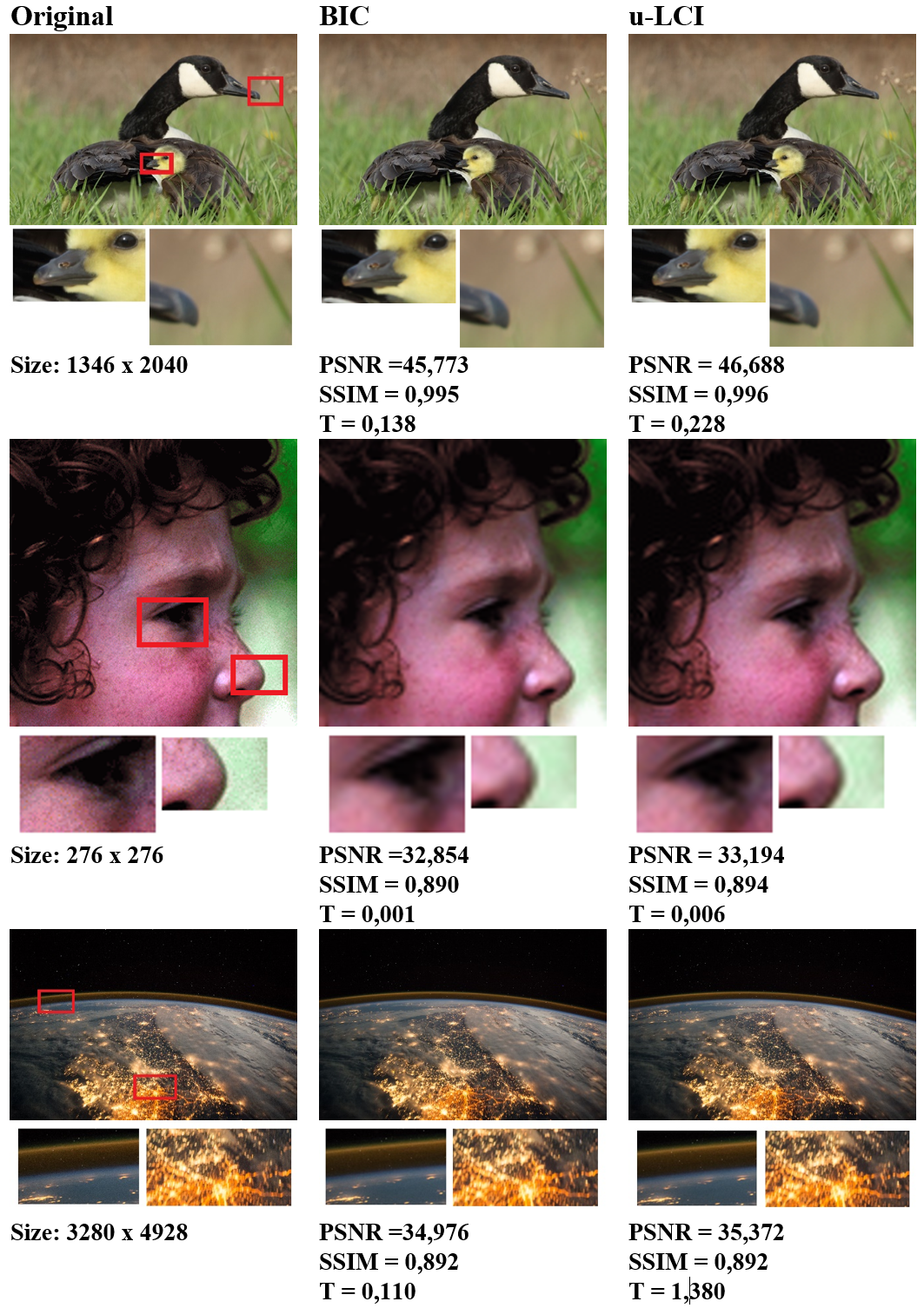}

\caption{Examples of upscaling performance results at the scale factor 2 (top),  at the scale factor 3 (middle), at the scale factor 4 (bottom).}
\label{fig:7}       
\end{center}
\end{figure}

\begin{figure}[h]
\begin{center}

\includegraphics[width=14cm]{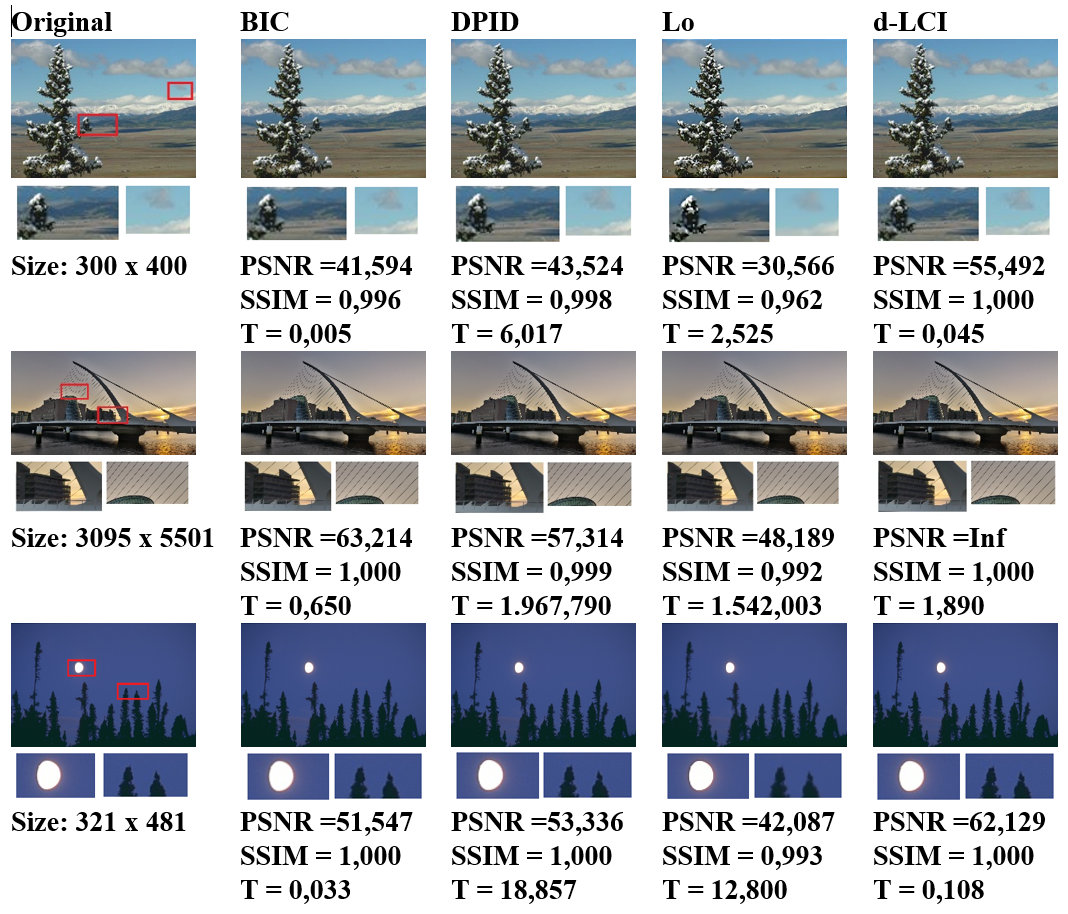}

\caption{Examples of downscaling performance results at the scale factor 2 (top),  at the scale factor 3 (middle), at the scale factor 4 (bottom).}
\label{fig_2}       
\end{center}
\end{figure}

As first experiment, we focus on the visual comparison of some original images chosen in the previous datasets with the output images produced by the methods in Subsection 5.1. Such images are given in Figure \ref{fig:7} for upscaling and in Figure \ref{fig_2} for downscaling.

In both the figures at the first column we show the target images and at the successive columns we display the resized images obtained as the output of the considered methods. Moreover, the first, second and third row of both the figures correspond to the scaling factors $2,3$ and $4$, respectively. In all the cases the images have been reported with evidence of some  magnified regions of interest (ROI).

From the qualitative point of view, by inspecting  Figures \ref{fig:7} and \ref{fig_2}, we  can deduce that  in terms of visual quality  LCI  has a good performance, also with respect to the chosen comparing methods.  It preserves the visual structure of the object without losing image details, the local contrast, and the luminance of the input image, by generating  resized images close to the target ones. Ringing and over smoothing artifacts are  limited and the resized images seem  fair,   sufficiently not blurred, and
 with  well-balanced colors.
\subsection{\it Quantitative comparison: average results}
For each dataset selected in Subsection 5.2, we computed the averages of the PSNR and SSIM values achieved by LCI and by  the comparison methods described in Subsection 5.1, for the scaling factors  $2,3,4$. The results are displayed in Table \ref{tab:9} for upscaling and in Table \ref{tab:2}  for downscaling, reporting in both the tables also the  averages of the required CPU times.
\begin{table}[ht]
\caption{Average  results in upscaling.}
\label{tab:9}       
\scriptsize{
\centering{\begin{tabular}{r|rrr|rrr|rrr}
\hline
  &\multicolumn{3}{c}{x2} & \multicolumn{3}{|c}{x3} & \multicolumn{3}{|c}{x4}\\ \hline
 & PSNR & SSIM & T  & PSNR & SSIM & T  & PSNR & SSIM & T\\
\hline
\textbf {BSDS500} & & & & & & & &\\
BIC  &	26,341&	0,886	&\textbf{0,003} 	&24,821	&0,837	&\textbf{0,002}&22,251	&0,701	 &\textbf{0,003} \\
u-LCI 	&\textbf{26,381}	&\textbf{0,888}	&0,014 	&\textbf{24,868}	&\textbf{0,839}	&0,010 &\textbf{22,446}&\textbf{0,769}	 &0,008\\
\hline
\textbf {NASA} & & & & & & & &\\
BIC &34,806 &	0,958	&\textbf{0,091} &	31,205 &	0,924	&\textbf{0,074}&29,808&	0,907 &\textbf{0,071} \\
u-LCI &\textbf{35,412} &	\textbf{0,960}	&1,357&	\textbf{31,542} &	\textbf{0,924	 }&0,839&\textbf{30,183}&	 \textbf{0,908} &0,634 \\
\hline
\textbf {YAHOO} & & & & & & & &\\
BIC  &	33,391	&	0,953&\textbf{0,056}&	29,754 &	0,913	&\textbf{0,055 }&28,741&	 0,891 &\textbf{0,051} \\
u-LCI 	&\textbf{33,900} &	\textbf{0,955	}&0,812  &\textbf{29,956} &\textbf{0,914}&0,567&\textbf{28,986}&	 \textbf{0,891} &0,459 \\
\hline
\textbf {URBAN100}& & & & & & & & \\
BIC & 25,433 &	0,882 &	\textbf{0,009} & 21,306 &	0,755 &	\textbf{0,008} & 21,703 &	0,741 &	 \textbf{0,008}\\
u-LCI & \textbf{25,896}	  & \textbf{0,886} &	0,064&		\textbf{21,325} &	\textbf{0,754} &	 0,051 & \textbf{21,919} &	\textbf{0,793} &	0,059
\\
\hline
\textbf {13US} & & & & & & & &\\
BIC 	 &	24,116 &	0,861 &	\textbf{0,002} & 20,754 &	0,734	&\textbf{0,002} &20,545 &	 0,710&	 \textbf{0,002}\\
u-LCI 	 &  \textbf{24,486} &	\textbf{0,868} &	0,012 & \textbf{20,780} &	\textbf{0,738} &	 0,010 & \textbf{20,656} &	\textbf{0,713} &	0,009\\
\hline
\end{tabular}
}}
\end{table}
\begin{table}[h]
\caption{Average results in downscaling }
\label{tab:2}
\scriptsize{
\centering{\begin{tabular}{r|rrr|rrr|rrr}
\hline
&\multicolumn{3}{c}{:2} & \multicolumn{3}{|c}{:3} & \multicolumn{3}{|c}{:4}\\ \hline
& PSNR & SSIM & T & PSNR & SSIM & T &  PSNR & SSIM & T\\
\hline
\textbf {BSDS500} & & & & & & & &\\ \hline
BIC     &38,872&	0,993	&\textbf{0,006} 	&39,253	&0,993	&\textbf{0,009 }&39,183	&0,993	 &\textbf{0,017} \\
DPID	&41,745&	0,996   &7,696  &41,827	&0,996	&12,246 	&41,206	&0,996	&18,615 \\
L$_0$	    &29,317&    0,961	&3,647 	&32,873	&0,971	&8,020 &34,174	&0,971	&14,295\\
d-LCI &\textbf{53,732}	&\textbf{1,000}	&0,057 	&\textbf{Inf}	&\textbf{1,000}	&0,091 	 &\textbf{55,889	 }&\textbf{1,000}	 &0,137\\
\hline
\textbf {NASA} & & & & & & & &\\
BIC &45,969	&0,995&	\textbf{0,229}&	47,114&	0,996&	\textbf{0,325}&	46,973&	0,996&	 \textbf{0,639}\\
DPID&	47,675	&0,998&	448,023&	47,657&	0,998&	731,498&	47,112&	0,997&	1.098,113\\
L$_0$ 	& 34,754& 0,972	& 208,574 &	37,285&	0,979&	617,386  & oom & oom & oom \\
d-LCI  	&\textbf{54,265}&	\textbf{0,999}&	6,614 	&\textbf{Inf}&	\textbf{1,000}&	 13,317&\textbf{	 57,491}&	 \textbf{1,000}&	 22,859\\
\hline
\textbf {YAHOO} & & & & & & & &\\
BIC &	44,757&	0,996	&\textbf{0,155} &	45,858	&0,997&	\textbf{0,219} 	&45,682&	0,996&	 \textbf{0,422}\\
DPID&	46,743&	0,998&	291,685 &	46,948	&0,998&	479,638&	46,421&	0,998&	714,908\\
L$_0$  &	33,913&	0,974&	133,190  & oom & oom & oom &  oom  & oom & oom\\
d-LCI 	&\textbf{54,067}&	\textbf{0,999}&	4,698  &	\textbf{Inf}	&\textbf{1,000}&	 8,279 	 &\textbf{57,223}&	 \textbf{1,000}&	 13,898\\
\hline
\textbf {URBAN00} & & & & & & & &\\
BIC &	35,661&	0,989	&\textbf{0,027} &	36,010	&0,990&	\textbf{0,041} 	&35,940&	0,989&	 \textbf{0,068}\\
DPID&	39,178&	0,996&	37,592 &	39,178	&0,996&	60,834&	38,702&	0,995&	93,996\\
L$_0$  &	26,718&	0,951&	13,267  & 31,061&	0,969 &	28,845 & 33,571 &	0,973 &	50,312\\
d-LCI 	&\textbf{52,613}&	\textbf{0,999}&	0,234  &	\textbf{Inf}	&\textbf{1,000}&	 0,416 	 &\textbf{55,452}&	 \textbf{1,000}&	 0,618\\
\hline
\textbf {13US} & & & & & & & &\\
BIC 	&35,129 &	0,990	&\textbf{0,005} &	35,469 &	0,991	&\textbf{0,009 }	 &35,397&	 0,991 &\textbf{	 0,013} \\
DPID	&38,061 &	0,996 &	5,593 	&38,326 &	0,996	&8,905 &37,707	&0,995	&13,819 \\
L$_0$ 	    &25,521	&0,949 &	2,572  &	30,231	&0,972	&5,706 &32,929	&0,979	&9,939\\
d-LCI   &\textbf{52,813} &	\textbf{1,000} &	0,042  &	\textbf{Inf}	&\textbf{1,000}	 &0,073  &\textbf{55,374	 }&\textbf{1,000}	&0,108 \\
\end{tabular}
}}
\end{table}

As announced in the introduction, from Table \ref{tab:9} we see that the
values attained by u-LCI are only slightly better than those achieved by BIC, and the execution time of u-LCI is a little bit greater than BIC.
Of course, as the scale factor increases,  both PSNR and SSIM decrease in all the upscaling methods. Moreover, the gaps between BIC and u-LCI  in both the metrics are maintained also for large size images.

 On the contrary, in downscaling, by inspecting Table \ref{tab:2} we observe that  the performance provided by d-LCI is much better than those achieved by the other methods for both the  metrics.
In particular, when the downscaling factor is $3$ 
 the theoretical results claimed in Proposition \ref{prop} are confirmed. Moreover, looking at the even downscaling factors $2,4$ 
 the improvement by d-LCI seems to increase as the scale factor increases while BIC seems to be affect by saturation, with an increasing gap  between d-LCI and the remaining methods.

 About the CPU time,  except for BIC,  d-LCI is  faster than the other methods that, in case of very large images, need very long computational time (see the results for NASA and YAHOO datasets that include target images of size $6394\times 3456$, i.e. input images with size $25576\times 13824$ in case of downscaling factor equals to 4) and in some cases, wherever we see oom (which means out of memory),  the available code of L$_0$ does not arrive to produce any result.
 \subsection{\it  Quantitative comparison: some pointwise results}
\begin{figure} [h]
\begin{center}

\includegraphics [width=15cm]{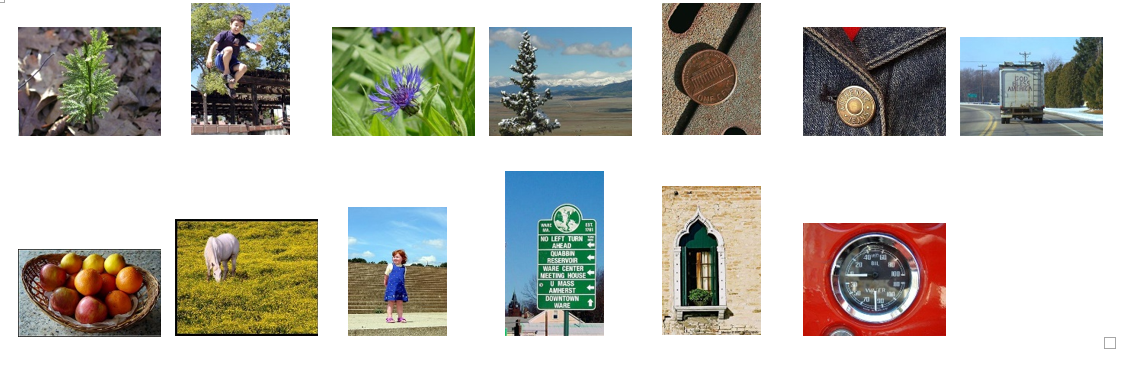}
\caption{13US dataset image: [U1 -U13] from left to right and from top to bottom
.}
\label{fig:1}       
\end{center}
\end{figure}

\begin{figure}[h]
\begin{center}

\includegraphics[width=15 cm]{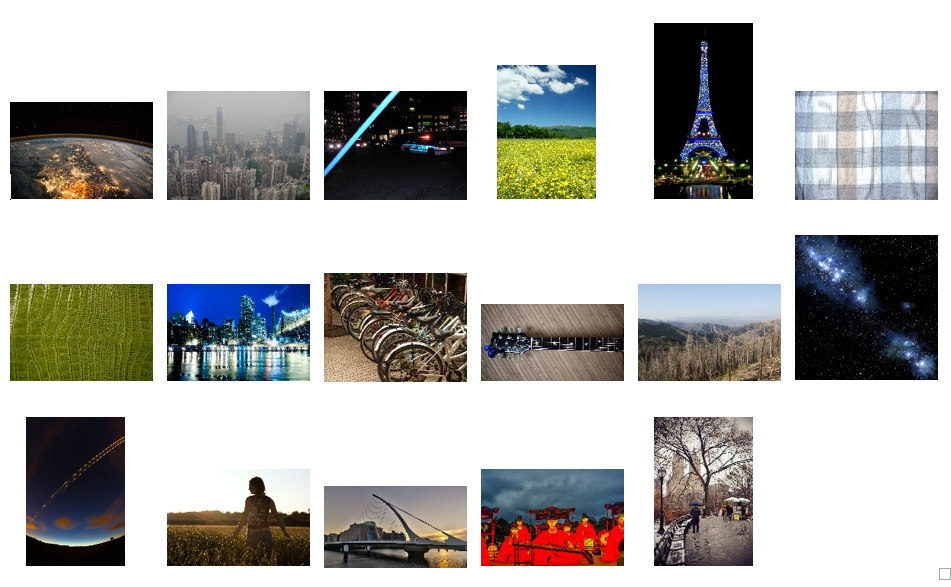}

\caption{NASA dataset images: [N1, N2, N3, N5, N6, N7, N8, N9, N10, N11, N12, N13, N14, N15, N16, N17, N19] from left to right and from top to bottom
.}
\label{fig:2}       
\end{center}
\end{figure}

 Here we focus on each single image from the smaller datasets of Table \ref{tab:1} , namely we consider 13US dataset in downscaling, with target images of small size displayed in Figure \ref{fig:1}, and NASA dataset in upscaling, with target images of large size displayed in Figure \ref{fig:2}.

The PSNR values achieved for every single image have been plotted in Figure 7 for NASA dataset with upscaling factor $2,4,8,16$, and in Figure 8 for 13US dataset with downscaling factor $6,18,30$.
\begin{figure}[!htbp]
\label{barr}
\begin{minipage}{.4\textwidth}
\includegraphics[scale=0.30]{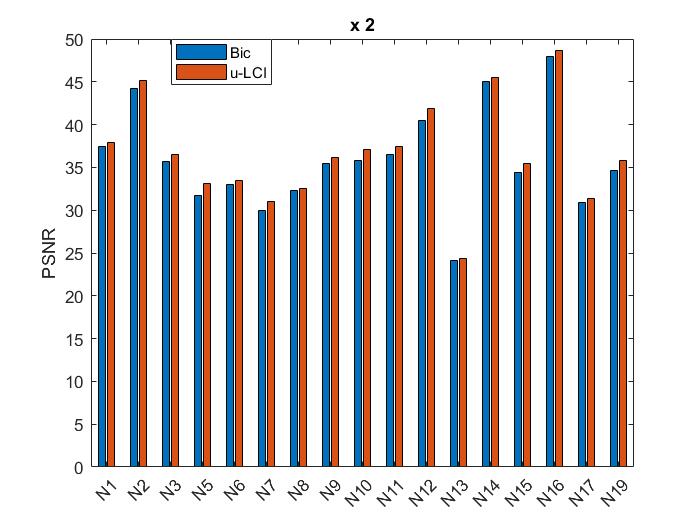}
\includegraphics[scale=0.30]{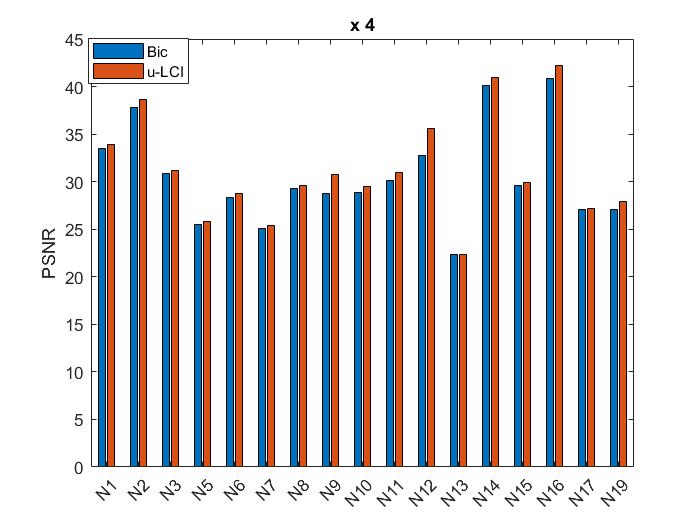}
\end{minipage}
\hspace{0.1cm}
\begin{minipage}{.4\textwidth}
\includegraphics[scale=0.30]{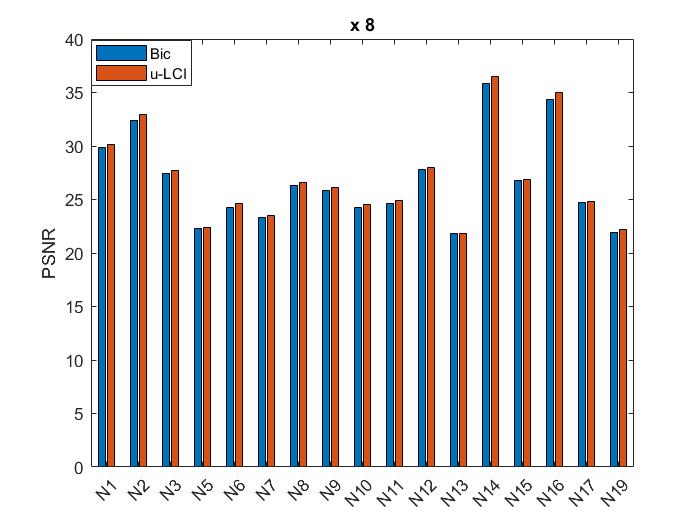}
\includegraphics[scale=0.30]{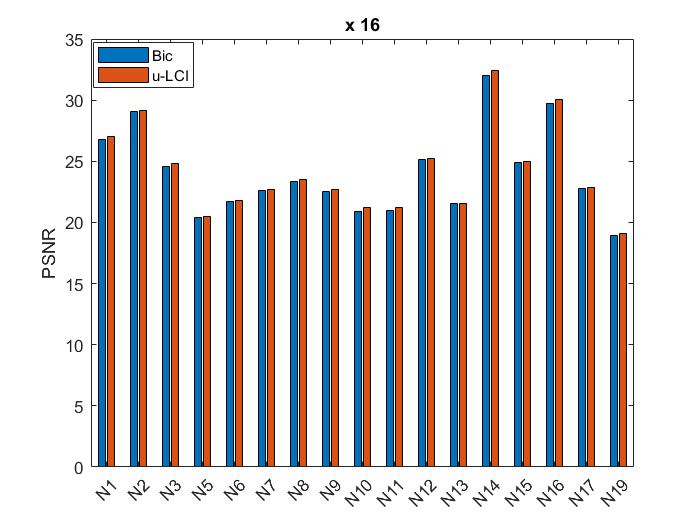}
\end{minipage}
\caption{Pointwise PSNR for the NASA set for scale factors   s=2,4,8,16}
\end{figure}
\begin{figure}[!htbp]
\label{stemm}
\begin{center}
\includegraphics[scale=0.28]{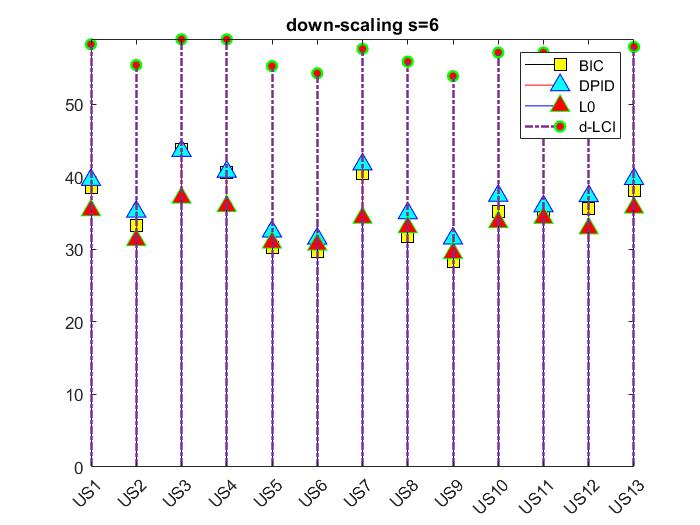}
\includegraphics[scale=0.28]{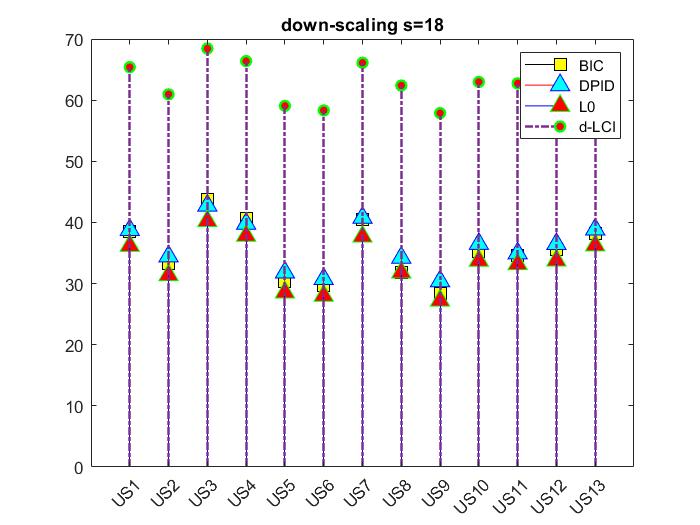}
\includegraphics[scale=0.28]{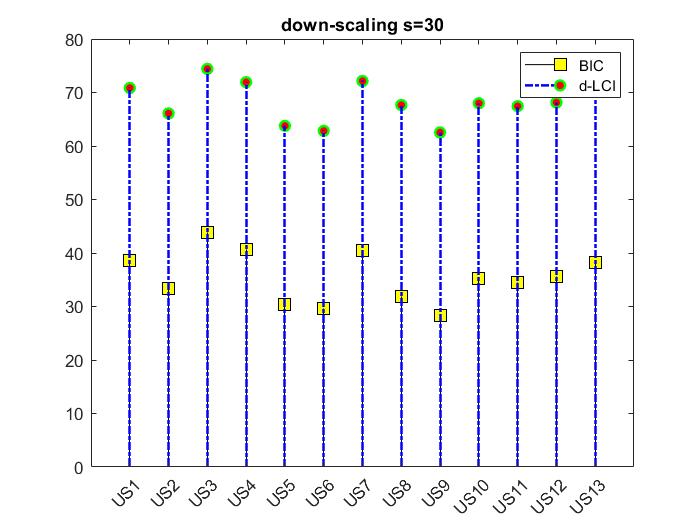}
\caption{Pointwise PSNR values for the set 13US, with scale factors s=6 (top), s=18 (middle), s=30 (down)}
\end{center}
\end{figure}

On the same datasets, more detailed results are given in Tables \ref{tab:11}--\ref{tab:5} reporting in the first columns the name and the size of all the images from 13US  and NASA.
\begin{table}[!htbp]
\caption{Pointwise results on NASA dataset  in upscaling}
\label{tab:11}       
\scriptsize{
\centering{\begin{tabular}{r|rr|rr|rr|rr}
\hline
  &\multicolumn{2}{c}{x2}  & \multicolumn{2}{|c}{x4} &\multicolumn{2}{c}{x8} &\multicolumn{2}{c}{x16}\\ \hline
& PSNR & SSIM   & PSNR & SSIM & PSNR & SSIM & PSNR & SSIM\\
\hline
\textbf {N1} \ 4928$\times$3280 & & & & & & \\
BIC& 37,472  &    0,945  & 33,453 & 0,892 &	29,828&	\textbf{0,858}		&26,795&	 \textbf{0,833}\\
u-LCI&\textbf{37,903}&\textbf{0,950}&\textbf{33,856}&0,892 &\textbf{30,133}&0,855&\textbf{26,999}&0,831\\
\hline
\textbf {N2} \ 3072$\times$2304 & & & & && \\
BIC&44,235&0,986&37,820&0,952 &32,417&0,888&29,065&0,840\\
u-LCI&\textbf{45,203}&\textbf{0,988}&\textbf{38,621}&\textbf{0,957} &\textbf{32,943}&\textbf{0,894}&\textbf{29,183}&0,840\\
\hline
\textbf {N3}\ 2048$\times$1536 & & & & && \\
BIC&35,708&0,963&30,845&\textbf{0,919}&27,474&\textbf{0,874}&24,593&\textbf{0,830}\\
u-LCI&\textbf{36,552}&\textbf{0,965}&\textbf{31,136}&0,916&\textbf{27,731}&0,865&\textbf{24,850}&0,808\\
\hline
\textbf {N5}\ 2701$\times$3665 & & & & &&\\
BIC&31,743&0,987&{25,444}&{0,960}&{22,262}&{0,932}&{20,388}&{0,913}\\
u-LCI&\textbf{33,077}&\textbf{0,989}&\textbf{25,792}&\textbf{0,962} &\textbf{22,426}&\textbf{0,933}&\textbf{20,487}&0,913\\
\hline
\textbf {N6} \ 1836$\times$3264 & & & & && \\
BIC&33,011&0,986&28,291&\textbf{0,965} &24,250&\textbf{0,927}&21,677&\textbf{0,886}\\
u-LCI&\textbf{33,443}&0,986&\textbf{28,779}&0,961&\textbf{24,614}&0,908&\textbf{21,784}&0,850\\
\hline
\textbf {N7}\ 1600$\times$1200  & & & & &&  \\
BIC&29,939&0,886&{25,087}&0,670&23,364&0,567&22,627&0,545\\
u-LCI&\textbf{30,981}&\textbf{0,905}&\textbf{25,409}&\textbf{0,684}&\textbf{23,542}&\textbf{0,569}&\textbf{22,683}&\textbf{ 0,546}\\
\hline
\textbf {N8}\ 3008$\times$2000 & & & &  &&\\
BIC&32,298&0,987&29,295&0,974  &26,346&0,956&23,345&0,932\\
u-LCI&\textbf{32,577}&\textbf{0,988}&\textbf{29,594}&\textbf{0,976}&\textbf{26,619}&\textbf{0,958}&\textbf{23,508}&\textbf{0,933}\\
\hline
\textbf {N9} \ 5430$\times$3520 & & & & && \\
BIC&35,440&0,985&28,707&0.945 &24,336&{0,883}&{22,494}&\textbf{0,846}\\
u-LCI&\textbf{36,116}&\textbf{0,986}&\textbf{30,735}&\textbf{0,955} &\textbf{26,121}&\textbf{0,899}&\textbf{22,712}&0,845\\
\hline
\textbf {N10} \ 3264$\times$2448 & & & & && \\
BIC&35,758&0,973&28,875&0,916 &24,213&0,837&20,914&0,767\\
u-LCI&\textbf{37,047}&\textbf{0,977}&\textbf{29,531}&\textbf{0,920}&\textbf{24,529}&\textbf{0,838}&\textbf{21,187}&\textbf{0,768}\\
\hline
\textbf {N11} \ 4368$\times$2326  & & & &  &&\\
BIC&36,525&\textbf{0,982}&30,149&\textbf{0,955}&24,649&\textbf{0,906}&20,966&\textbf{0,849}\\
u-LCI&\textbf{37,461}&0,981&\textbf{30,986}&{0,951}&\textbf{24,918}&0,897&\textbf{21,193}&0,846\\
\hline
\textbf {N12}\ 3504$\times$2336 & & & &  &&\\
BIC&40,471&0,989&32,786&0,951 &27,777&0,888&25,116&0,850\\
u-LCI&\textbf{41,889}&\textbf{0,991}&\textbf{35,581}&\textbf{0,956}&\textbf{27,982}&\textbf{0,890}&\textbf{25,191}&\textbf{0,851}\\
\hline
\textbf {N13}\ 1200$\times$1600 & & & & && \\
BIC&24,145&\textbf{0,746}&22,370&\textbf{0,603} &21,838&\textbf{0,530}&21,555&0,483\\
u-LCI&\textbf{24,316}&0,730&\textbf{22,390}&0,592  &\textbf{21,852}&0,529&\textbf{21,579}&\textbf{0,486}\\
\hline
\textbf {N14} \ 3456$\times$5184 & & & & &&  \\
BIC&45,041&0,992&40,137&0,986  &35,808&0,978&31,982&\textbf{0,967}\\
u-LCI&\textbf{45,539}&0,992&\textbf{40,907}&0,986 &\textbf{36,531}&0,978&\textbf{32,407}&0,966\\
\hline
\textbf {N15} \ 3072$\times$2048 & & & &  &&\\
BIC&34,457&0,990&29,574&0,967  &26,731&\textbf{0,941}&24,886&\textbf{0,920}\\
u-LCI&\textbf{35,427}&\textbf{0,991}&\textbf{29,872}&\textbf{0,968}&\textbf{26,883}&0,940&\textbf{24,976}&0,919\\
\hline
\textbf {N16}\ 5501$\times$3095  & & & & && \\
BIC&47,929  &    0,993&40,846&0,982  &34,304&0,961&29,760&\textbf{0,940}\\
u-LCI &\textbf{48,691}&\textbf{0,994}&\textbf{42,202}&\textbf{0,983}&\textbf{35,029}&{0,961}&\textbf{30,027}&0,939\\
\hline
\textbf {N17} \ 2048$\times$1363 & & & & && \\
BIC&30,892 &     0,980&27,028&{0,954} &24,688&{0,927}&\textbf{22,895}&\textbf{0,899}\\
u-LCI&\textbf{31,380}&\textbf{0,981}&\textbf{27,219}&{0,954}&\textbf{24,834}&0,927&21,044&0,877\\
\hline
\textbf {N19} \ 3039$\times$4559 & & & & && \\
BIC&34,581   &   0,977&27,048&0,884  &21,873&0,718&18,910&0,615\\
u-LCI&\textbf{35,844}&\textbf{0,981}&\textbf{27,881}&\textbf{0,896}&\textbf{22,178}&\textbf{0,723}&\textbf{19,077}&\textbf{0,617}\\
\hline
\end{tabular}
}}
\end{table}
\begin{table}[!htbp]
\caption{Pointwise results on 13US dataset  in downscaling}
\label{tab:3}       
\scriptsize{
\centering{\begin{tabular}{r|rr|rr|rr|rr}\hline
  &\multicolumn{2}{c}{:3} & \multicolumn{2}{|c}{:6} & \multicolumn{2}{|c}{:9} & \multicolumn{2}{|c}{:18}
  \\ \hline
 & PSNR & SSIM  & PSNR & SSIM   & PSNR & SSIM & PSNR & SSIM \\
\hline
\textbf {US1}\  300$\times$400 & &     & & & & & &
\\
BIC &38,625&0,996 &38,573 &	0,996 & 38,583 &	0,996 & 38,584&	0,996\\
DPID&40,846&	0,997&39,580 &	0,997 & 39,160 &	0,996 & 38,764&	0,996\\
L$_0$ &27,994&0,949& 35,405	&0,983 & 35,138 &	0,986 & 36,183 &	0,991\\
d-LCI &\bf{Inf} &\bf{1,000}  &	\textbf{58,265} &	\textbf{1,000} 	&\bf{Inf} &\bf{1,000} &	 \textbf{65,394} &	 \textbf{1,000}\\
\hline
\textbf {US2} \  400$\times$300 & &  & & & & & &
\\
BIC&33,348&0,987& 33,298	&0,987 &33,303&	0,987 &33,302	&0,987\\
DPID&36,428&0,994 & 35,217 &	0,992 & 34,823&	0,992&  34,450&	0,991\\
L$_0$ &29,184&0,963& 31,289&	0,967 & 30,647	&0,968& 31,394	&0,977\\
d-LCI&\bf{Inf}&\bf{1,000}&\textbf{55,418}	&\textbf{1,000} &\bf{Inf}&\bf{1,000} &\textbf{60,956}&	 \textbf{1,000}\\
\hline
\textbf {US3} \  300$\times$400 & &  & & & & & &
\\
BIC&43,820&0,999   &43,754&	0,999& 43,774&	0,999& 43,786&	0,999 \\
DPID&44,853&0,999 & 43,589&	0,999&43,162&	0,999& 42,751	&0,999\\
L$_0$ &34,756&0,992 &  37,140&	0,995&38,087	&0,996&40,237&	0,998\\
d-LCI&\bf{Inf}&\bf{1,000}&\textbf{58,958}&	\textbf{1,000 }& \textbf{Inf}	& \textbf{1,000 }&\textbf{ 68,436}	 &\textbf{1,000}
\\
\hline
\textbf {US4} \  300$\times$400 & &  & & & & & &
\\
BIC&40,673&0,996 &40,623	&0,996 & 40,631 &	0,996 & 40,633&	0,996\\
DPID&42,066&0,998 &40,699&	0,997& 40,236	&0,997 & 39,794	&0,996\\
L$_0$ &32,944&0,973 &35,973&	0,981& 36,453	&0,986& 37,836	&0,992\\
d-LCI&\textbf{Inf}&\textbf{1,000}& \textbf{58,958}&	\textbf{1,000}& \textbf{Inf}	 &\textbf{1,000}&\textbf{ 66,386}&	 \textbf{1,000}\\
\hline
\textbf {US5} \  400$\times$300 & &  & & & & & &
\\
BIC&30,357&0,979 & 30,307	&0,979& 30,311	& 0,979& 30,311	&0,979\\
DPID&33,492&0,991& 32,471&	0,989& 32,133	&0,988& 31,813	&0,987\\
L$_0$ &26,807&0,962& 30,944	&0,973&28,911&	0,963&28,551&	0,965\\
d-LCI &   \textbf{  Inf}&   \textbf{ 1,000}& \textbf{ 55,257}&	\textbf{1,000}&\textbf{Inf}	 &\textbf{1,000}&\textbf{59,065}&	\textbf{1,000}
\\
\hline
\textbf {US6} \  300$\times$400 & &  & & & & & &
\\
BIC&29,718&0,978 & 29,659&	0,978&29,633&	0,978&29,663	&0,978\\
DPID&32,578&0,990& 31,486&	0,988&31,111&	0,987 & 30,752	&0,985\\
L$_0$ &25,322&0,959&30,666	& 0,966 &28,724	& 0,958 &28,022	&0,962\\
d-LCI&\textbf{Inf}&  \textbf{  1,000}& \textbf{54,273}	& \textbf{1,000 }& \textbf{Inf}	 &\textbf{1,000} &\textbf{ 58,321}	&\textbf{1,000}\\
\hline
\textbf {US7} \  273$\times$400  & &  & & & & & &
\\
BIC&40,559&0,995 & 40,513 &	0,995 &40,521	&0,995 &40,522&	0,995 \\
DPID &43,135&0,997& 41,717	&0,996& 41,232	&0,996 & 40,771& 	0,996\\
L$_0$ &31,832&0,949& 34,364	&0,964 & 35,604&	0,975& 37,697&	0,987\\
d-LCI& \textbf{  Inf}&  \textbf{ 1,000}& \textbf{57,639}& 	\textbf{1,000 }& \textbf{Inf} &	 \textbf{1,000} &\textbf{66,113}	&\textbf{1,000}\\
\hline
\textbf {US8}  \  322$\times$400& &  & & & & & &
\\
BIC&31,857&0,989 & 31,802 &	0,989 &31,808 &	0,989 &31,808	&0,989 \\
DPID&36,382&0,996 & 34,979&	0,994&34,601&	0,994&34,239	&0,993\\
L$_0$&29,562&0,974 & 33,051	&0,977 &32,378	&0,980 &31,828 &	0,985\\
d-LCI&  \textbf{  Inf}&    \textbf{1,000} & \textbf{55,880}	& \textbf{1,000} &\textbf{Inf}	& \textbf{1,000} &\textbf{62,394} &	\textbf{1,000}\\
\hline
\textbf {US9} \  322$\times$400 & &  & & & & & &
\\
BIC&28,407&0,984 & 28,350&	0,984 &28,536 &	0,984 & 28,354& 0,984\\
DPID&32,542&0,994& 31,504&	0,993 & 34,601 &	0,992 & 30,393&	0,992\\
L$_0$&29,562&0,974 &29,504	&0,985 & 27,361 &	0,980 & 27,214& 	0,979\\
d-LCI&  \textbf{ Inf}&    \textbf{1,000}& \textbf{53,879 }&	\textbf{1,000} & \textbf{Inf}&\textbf{1,000} & \textbf{57,877}&	\textbf{1,000}\\
\hline
\textbf {US10} \  400$\times$307 & &  & & & & & &
\\
BIC&35,220&0,992 & 35,168	&0,992&35,172	&0,992&35,173&	0,992\\
DPID&38,947&0,997 & 37,406&	0,996 &36,929	& 0,995 &36,494	& 0,995\\
L$_0$ &31,031&0,984  & 33,760	& 0,986   &33,227 &	0,985 & 33,744&	0,988\\
d-LCI&\textbf{Inf}&\textbf{1,000} & \textbf{57,131}&	\textbf{1,000}& \textbf{Inf	}&\textbf{ 1,000} & \textbf{62,986	 }& \textbf{1,000}
\\
\hline
\textbf {US11} \  400$\times$241 & &  & & & & & &
\\
BIC& 34,599&0,995 & 34,528	&0,995 & 34,536	& 0,995 &34,535	& 0,995\\
DPID&37,345&0,998 & 35,950 &	0,997 & 35,432 & 	0,997 &  34,937	& 0,996 \\
L$_0$&27,631&0,976 & 34,346 & 	0,986 & 34,143 & 	0,989 &  33,189 & 	0,992\\
d-LCI&\textbf{Inf}&\textbf{1,000}& \textbf{57,147}	& \textbf{1,000 }&\textbf{ Inf	} & \textbf{ 1,000} &  \textbf{62,744}	& \textbf{61,000}\\
\hline
\textbf {US12} \ 400$\times$266 & & & & & & & &
\\
BIC&35,698&0,993 & 35,639	&0,993 &035,645&	0,993 &35,647&	0,993\\
DPID& 38,626&0,997& 37,345&	0,996 &36,898	& 0,995 & 36,471	&0,995\\
L$_0$ &29,883&0,972& 32,902&	0,978 & 33,091&	0,982& 33,776&	0,987\\
d-LCI&\textbf{Inf}&\textbf{1,000}&\textbf{ 55,801} &	\textbf{1,000} & \textbf{Inf}	& \textbf{1,000} & \textbf{62,744}	& \textbf{1,000}\\
\hline
\textbf {US13} \ 310$\times$400 & &  & & & & & &
\\
BIC &38,220 &0,996  & 38,166 &	0,996 & 38,173& 	0,996 & 38,177	& 0,996\\
DPID&41,001&0,998& 39,707 &	0,997 &39,276 & 	0,997 &38,861	& 0,997 \\
L$_0$&31,593&0,982& 35,761& 	0,989 &35,898&	0,991 &36,263 &	0,993 \\
d-LCI&\textbf{Inf}&\textbf{1,000}& \textbf{57,905}	& \textbf{1,000} & \textbf{Inf} & \textbf{	 1,000} & \textbf{65,535}	& \textbf{1,000}
\\
\hline
\end{tabular}
}}
\end{table}
\begin{table}[!htbp]
\caption{Performance results on 13US dataset  for high downscaling}
\label{tab:5}       
\scriptsize{
\centering{\begin{tabular}{r|rr|rr|rr|rr|rr}
\hline
& \multicolumn{2}{|c}{:20} & \multicolumn{2}{|c}{:21} & \multicolumn{2}{|c}{:27}& \multicolumn{2}{|c}{:30}& \multicolumn{2}{|c}{:33}\\ \hline
 &  PSNR & SSIM   & PSNR & SSIM & PSNR & SSIM & PSNR & SSIM & PSNR & SSIM\\
\hline
\textbf {US1}\ 300$\times$400 & & & & &&&&\\
BIC &		38,583	&0,996  &	38,585&	0,996		&	38,586	&0,996              &	 38,585&	 0,996 &	38,584&	 0,996\\
d-LCI &		\textbf{66,485}&\textbf{1,000}     &	\textbf{Inf}	&\textbf{1,000} 	 &\textbf{Inf}&	\textbf{1,000	 }              &\textbf{70,845}&	 \textbf{1,000} &\textbf{Inf}&	 \textbf{1,000}\\
\hline
\textbf {US2} \ 400$\times$300 & & & & &&&&\\
BIC &		33,302 &	0,987	&	33,302 &	0,987  &	33,302 &	0,987	                       &	33,303&	 0,987&	33,303&	0,987\\
d-LCI&	\textbf{66,347}&	\textbf{1,000}   &	\textbf{Inf}	&\textbf{1,000}&	 \textbf{Inf}	                &\textbf{66,092}	&\textbf{1,000} &\textbf{1,000}&\textbf{Inf}	 &\textbf{1,000}\\
\hline
\textbf {US3} \ 300$\times$400 & & & & &&&&\\
BIC &43,786&0,999  &43,787&0,999&43,786&0,999&43,788&0,999&43,784&0,999\\
d-LCI&\textbf{69,721}&\textbf{1,000}&\textbf{Inf}	&\textbf{1,000}&\textbf{Inf}	 &\textbf{1,000} &\textbf{74,415}&\textbf{1,000}&\textbf{Inf}	&\textbf{1,000}\\
\hline
\textbf {US4} \ 300$\times$400 & & & & &&&&\\
BIC &	40,633&	0,996  &40,634&	0,996	&	40,636&	0,996                              &	 49,636&	0,996&	49,635&	 0,996\\
d-LCI&	\textbf{68,054}&	\textbf{1,000}       &	\textbf{Inf}	&\textbf{1,000}&	 \textbf{Inf}	&\textbf{1,000}        &	\textbf{71,936}&	 \textbf{1,000}&	 \textbf{Inf}	 &\textbf{1,000}\\
\hline
\textbf {US5} \ 400$\times$300 & & & & &&&&\\
BIC &30,311&0,979 &30,311&0,979&30,311&0,979                             &30,311&0,979&30,311&0,979\\
d-LCI&\textbf{59,853}&\textbf{1,000} &\textbf{Inf}&\textbf{1,000}&\textbf{Inf}&\textbf{1,000}   &\textbf{63,768}&\textbf{1,000}&\textbf{Inf}&\textbf{1,000}\\
\hline
\textbf {US6} \ 300$\times$400 & & & & &&&&\\
BIC &30,311&0,979  &29,663&0,978&29,664&0,978                              &29,663&0,978&29,663&0,978\\
d-LCI&\textbf{59,135}&\textbf{1,000} &\textbf{Inf}&\textbf{1,000}&\textbf{Inf}&\textbf{1,000}                            &\textbf{62,833}&\textbf{1,000}&\textbf{Inf}&\textbf{1,000}\\
\hline
\textbf {US7} \ 273$\times$400 & & & & &&&&\\
BIC &40,524&0,995 &40,521&0,995&40,524&0,995                             &40,536&0,995&40,520&0,995\\
d-LCI&\textbf{67,318}&\textbf{1,000} &\textbf{Inf}&\textbf{1,000}&\textbf{Inf}&\textbf{1,000}                           &\textbf{72,151}&\textbf{1,000}&\textbf{Inf}&\textbf{1,000}\\
\hline
\textbf {US8}  \ 322$\times$400 & & & & &&&& \\
BIC &31,808&0,989 &31,807&0,989&31,807&0,989                         &31,807&0,995&31,806&0,989\\
d-LCI&\textbf{63,539}&\textbf{1,000} &\textbf{Inf}&\textbf{1,000}&\textbf{Inf}&\textbf{1,000}                          &\textbf{67,653}&\textbf{1,000}&\textbf{Inf}&\textbf{1,000}\\
\hline
\textbf {US9}  \ 322$\times$400 & & & & &&&&\\
BIC &28,354&0,984  &28,355&0,984&28,354&0,984                                         &28,355 &0,984&28,355 &0,984\\
d-LCI&\textbf{58,651}&\textbf{1,000}  &\textbf{Inf}&\textbf{1,000}&\textbf{Inf}&\textbf{1,000}                       &\textbf{62,536}&\textbf{1,000}&\textbf{Inf}&\textbf{1,000}\\
\hline
\textbf {US10} \ 400$\times$307 & & & & &&&&\\
BIC&	35,174&	0,992&	35,172	&0,992&	35,174&	0,992                                      &		 35,175&	0,992&	 35,175&	0,992\\
d-LCI&\textbf{63,924}&	\textbf{1,000}&\textbf{Inf}&\textbf{1,000}&\textbf{Inf}&\textbf{1,000}                &	 \textbf{67,994}&	\textbf{1,000}&\textbf{Inf}&\textbf{1,000}\\
\hline
\textbf {US11} \ 400$\times$241 & & & & &&&&\\
BIC &	34,535&	0,995 &	34,535&	0,995&	34,534&	0,995     &34,534&	0,995 &34,535&	0,995\\
d-LCI	&	\textbf{63,445}&	\textbf{1,000} &\textbf{Inf}&\textbf{1,000}&\textbf{Inf}&\textbf{1,000}                             &	 \textbf{67,442}&\textbf{1,000}&\textbf{Inf}&\textbf{1,000}\\
\hline
\textbf {US12} \ 400$\times$266 & & & & &&&&\\
BIC &	35,647&	0,993  &	35,646	&0,993&	35,645&	0,993    &	35,646	&0,993&	35,648	 &0,993\\
d-LCI	&\textbf{63,760}&	\textbf{1,000} &\textbf{Inf}&\textbf{1,000}&\textbf{Inf}&\textbf{1,000}                           &	 \textbf{67,127}	 &\textbf{1,000}&\textbf{Inf}&\textbf{1,000}\\	
\hline
\textbf {US13} \ 310$\times$400 & & & & &&&&\\
BIC 	&	38,175&	0,996 &38,175	&0,996&	38,174&	0,996                                   &	 38,175	&0,996	&	 38,175	&0,996	\\
d-LCI &	\textbf{66,347} &	\textbf{1,000}  &\textbf{Inf}&\textbf{1,000}&\textbf{Inf}&\textbf{1,000}                      &\textbf{70,960} &	 \textbf{1,000}&\textbf{Inf}&\textbf{1,000}\\
\hline
\end{tabular}
} 
}
\end{table}

In particular, Table \ref{tab:11}  contains the detailed upscaling results on the NASA dataset,  for the scaling factors $ s\in \{2, 4, 8, 16\}$ while Tables \ref{tab:3}-\ref{tab:5} concern the downscaling results obtained on  13US dataset  for the scaling factors $s\in \{3,  6, 9, 18\}$ (Table \ref{tab:3}) and for the larger factors $s\in \{20, 21, 27, 30,  33\}$ (Table \ref{tab:5}).

 According to the average results,  the pointwise ones corroborate the global trend.
In particular, Table \ref{tab:3} confirms that in downscaling d-LCI is the method allowing the largest PSNR and SSIM, and that the method requiring the least computation time is BIC, followed by d-LCI which is faster than DPID and L$_0$. Also, note that the scaling factors chosen in Table \ref{tab:3} are all multiple of three but, according to Proposition \ref{prop} we can see the differnce between the odd and even downscaling factors. Moreover, we can affirm that d-LCI works fine also for large downscaling factors as reported in Table \ref{tab:5} where the target images in 13US have been zoomed in up to 33 times, getting input images whose size goes up to 13200$\times$10230.
We point out that in Table \ref{tab:5} the results for DPID and L$_0$ methods are missing  because the  publicly available MatLab codes  don't work for so large scale factors. Moreover, if $s\in \{21,27,33\}$ then, in accordance with  Proposition \ref{prop},  PSNR and SSIM again reach the limit values  Infinite and 1, respectively, by d-LCI  while the same does not hold for BIC applied to the same input images.
\section{Conclusions}
In the context 
of interpolation methods for  image resizing, we present the Lagrange--Chebyshev Interpolation (LCI) method. It is based on a non standard mathematical modeling that leads to the application of such an optimal Lagrange interpolation process to globally approximate the image at a continuous scale.
One of the main advantages of LCI method is its high flexibility of  working in both the scaling directions, either  setting a scale factor or giving a particular final size.
Comparisons with other resizing procedures have been reported on 5 different common datasets made   of 726 images in total. The numerical experience in upscaling shows a performance
comparable with the bicubic interpolation method, while in downscaling a much better performance is obtained with respect to all the considered comparison methods. Moreover, in downscaling cases with odd scale factors, we estimate the Mean Square Error (MSE) of our procedure in terms of the initial errors present in the data and we prove that it is null in absence of noise or artifacts in the input image.
We wonder whether further improvements can be achieved employing wavelets technique or finer approximation polynomials. These will be the subjects of further investigations.
\vspace{1.5cm}\newline
\textbf{Code and supplementary materials}

The code and the supplementary materials are openly available (to appear at a GitHub link).
\bibliographystyle{plain}      
\bibliography{biblio}   
\end{document}